\documentclass[11pt,bezier]{article}
\usepackage{amsmath, amssymb, amsfonts, graphicx}

\newtheorem{prethm}{{\bf Theorem}}

\newtheorem{prepro}[prethm]{{\bf Proposition}}

\newtheorem{prelem}[prethm]{{\bf Lemma}}

\newtheorem{predeff}[prethm]{{\bf Definition}}

\newtheorem{precor}[prethm]{{\bf Corollary}}

\newtheorem{preconj}[prethm]{{\bf Conjecture}}

\newenvironment{conj}{\begin{preconj}\sl{\hspace{-0.5
               em}{\bf.}}}{\end{preconj}}

\newtheorem{preremark}[prethm]{{\bf Remark}}

\newtheorem{preexample}[prethm]{{\bf Example}}

\newtheorem{preproof}{{\bf\textsf{Proof.}}}

\title{Nontrivial nuciferous graphs exist}

\author{{\sc Ebrahim Ghorbani} \\[.3cm]
{\sl Department of Mathematics, K.N. Toosi University of Technology,}\\
{\sl P. O. Box 16315-1618, Tehran, Iran}\\
{\sl School of Mathematics, Institute for Research in Fundamental
Sciences (IPM),}\\
{\sl P.O. Box
19395-5746, Tehran, Iran }
\\[.3cm]
$\mathsf{e\_ghorbani@ipm.ir}$ }

\begin{document}
\maketitle

\vspace{5mm}

\begin{abstract}
A nuciferous graph is a simple graph with a non-singular $0$-$1$ adjacency
matrix $A$ such that all the diagonal entries of $A^{-1}$ are zero and all
the off-diagonal entries of $A^{-1}$ are non-zero. Sciriha {\it et al.}
conjectured that except $K_2$, no nuciferous graph exists. We disprove this
conjecture. Moreover, we conjecture that there infinitely many nuciferous
Cayley graphs.
\vspace{5mm}

\noindent {\bf AMS Mathematics Subject Classification\,(2010):}   05C50
\end{abstract}

\vspace{5mm}

\section{Introduction}

Let $G$ be a simple graph with non-singular $0$-$1$ adjacency matrix $A$. If
all the diagonal entries of $A^{-1}$ are zero and all the off-diagonal
entries of  $A^{-1}$ are non-zero, then $G$ is called a {\em nuciferous graph}.
 The concept of nuciferous graphs has arisen in the context of
the quantum mechanical theory of the conductivity of non-singular
carbon molecules in the Source and Sink Potential model \cite{f}. According to Sciriha \cite{s0}:
\begin{quote}
 ``In the graph-theoretical Source and Sink Potential model, a molecule is either an insulator or a conductor for electrons with
energy zero. Of particular interest are two classes of graphs with analogous vertex pairs, i.e., the same behavior for any two-
vertex connection. These are uniform-core (insulating for all two-vertex connections) and nuciferous graphs, which conduct for
all two-vertex connections. A graph $G$ in the first class reaches the minimum possible nullity when any two distinct connecting
vertices are deleted. In the second class, the nullity reaches one, the maximum possible, when any vertex is deleted.''
\end{quote}
To date, the only nuciferous graph known is $K_2$ which we call it the trivial one. In \cite{s} it was conjectured
that there are no others. We disprove this
conjecture.

We remark that in \cite{fa} weighted graphs that have an adjacency matrix with the
required structure in the inverse as in  nuciferous graphs were found.

\section{Nuciferous Cayley graphs}

Making use of the database of vertex-transitive graphs by Gordon Royle \cite{r} and an exhaustive computer search, we found several nuciferous graphs.
In fact among vertex-transitive graphs with at most 31 vertices, there are 21 nuciferous graphs.
All the 21 nuciferous graphs we found are Cayley graphs: 6 on 24, 3 on 28 and 12 on 30 vertices.
Recall that a Cayley graph Cay$(\Gamma,
S)$ for a given group $\Gamma$ and connection set $S \subset\Gamma$ is the
graph with vertex set  $\Gamma $ and with $u$ connected to $v$ if and only
if $vu^{-1} \in S$.
Table~\ref{table} shows the list of nuciferous Cayley graphs up to 31vertices according to their groups in which $C_n$ and $D_n$ denote cyclic and dihedral groups of order $n$, respectively, and ${\rm Sym}_k$ and ${\rm Alt}_k$ denote the symmetric and alternative groups on $k$ elements, respectively.
\begin{table}[!ht]
\centering
\begin{tabular}{clcc}
  \hline
  Order & Group & Degree&\# of nuciferous \\
  \hline
24 &	$D_{12}\times C_2$	&15& 1 \\
24	& ${\rm Alt}_4\times C_2$& 7&	$2$ \\
24	 &	${\rm Sym}_3\times C_4$ &15&	1\\
24		& $D_{24}$&	15&1 \\
24		&${\rm Sym}_4$	&7&2 \\
24		&${\rm Sym}_4$	&15&1\\\hline
28		&$D_{28}$&13&1\\
28		&$D_{28}$&15&2\\\hline
30	&	$C_{30}$&11& 1\\
30	&	$C_{30}$&15& 1\\
30		&$D_{10}\times C_3$&11&1\\
30		&$D_{10}\times C_3$&22&1\\
30		&$D_6\times C_5$	&11&1\\
30		&$D_{30}$&11&1\\
30		&$D_{30}$&15&10\\
\hline\end{tabular}
\caption{Number of nuciferous Cayley graphs with at most 31 vertices  according to their group}\label{table}
\end{table}
We notice that in Table~\ref{table} the two degree 7 graphs on ${\rm Alt}_4\times C_2$ are isomorphic to the two degree 7 on ${\rm Sym}_4$; all the degree 11 graphs are isomorphic; and the degree 15 graph on $C_{30}$ is isomorphic to one of the degree 15 graphs on  $D_{30}$.

  Table~\ref{table2} depicts the adjacency matrix $A$ and the inverse $A^{-1}$  of one of the two nuciferous Cayley graphs with 24 vertices on the group
  ${\rm Alt}_4\times C_2$.
\begin{table}[!ht]
$$A={\footnotesize\left[\begin{array}{cccccccccccccccccccccccc}
0& 1& 1& 1& 1& 1& 1& 1& 0&0 &0 &0 &0 &0 &0 &0 &0 &0 &0 &0 &0 &0 &0 & 0\\
1& 0& 1& 1& 0& 0& 0& 0& 0& 0& 0& 0& 1& 1& 0& 0& 0& 0& 0& 0& 1& 1& 0& 0\\
1& 1& 0& 1& 0& 0& 0& 0& 1& 0& 1& 0& 0& 0& 1& 0& 0& 0& 1& 0& 0& 0& 0& 0\\
1& 1& 1& 0& 0& 0& 0& 0& 0& 1& 0& 1& 0& 0& 0& 1& 0& 0& 0& 1& 0& 0& 0& 0\\
1& 0& 0& 0& 0& 1& 0& 0& 0& 1& 0& 1& 0& 0& 1& 0& 0& 1& 0& 0& 1& 0& 0& 0\\
1& 0& 0& 0& 1& 0& 0& 0& 1& 0& 1& 0& 0& 0& 0& 1& 1& 0& 0& 0& 0& 1& 0& 0\\
1& 0& 0& 0& 0& 0& 0& 0& 0& 1& 0& 1& 1& 1& 0& 0& 1& 0& 1& 0& 0& 0& 0& 0\\
1& 0& 0& 0& 0& 0& 0& 0& 1& 0& 1& 0& 1& 1& 0& 0& 0& 1& 0& 1& 0& 0& 0& 0\\
0& 0& 1& 0& 0& 1& 0& 1& 0& 1& 0& 0& 0& 1& 0& 0& 0& 0& 0& 1& 0& 0& 1& 0\\
0& 0& 0& 1& 1& 0& 1& 0& 1& 0& 0& 0& 1& 0& 0& 0& 0& 0& 1& 0& 0& 0& 1& 0\\
0& 0& 1& 0& 0& 1& 0& 1& 0& 0& 0& 0& 0& 0& 1& 1& 0& 0& 0& 0& 0& 1& 0& 1\\
0& 0& 0& 1& 1& 0& 1& 0& 0& 0& 0& 0& 0& 0& 1& 1& 0& 0& 0& 0& 1& 0& 0& 1\\
0& 1& 0& 0& 0& 0& 1& 1& 0& 1& 0& 0& 0& 0& 1& 0& 0& 1& 1& 0& 0& 0& 0& 0\\
0& 1& 0& 0& 0& 0& 1& 1& 1& 0& 0& 0& 0& 0& 0& 1& 1& 0& 0& 1& 0& 0& 0& 0\\
0& 0& 1& 0& 1& 0& 0& 0& 0& 0& 1& 1& 1& 0& 0& 0& 0& 0& 0& 0& 1& 0& 1& 0\\
0& 0& 0& 1& 0& 1& 0& 0& 0& 0& 1& 1& 0& 1& 0& 0& 0& 0& 0& 0& 0& 1& 1& 0\\
0& 0& 0& 0& 0& 1& 1& 0& 0& 0& 0& 0& 0& 1& 0& 0& 0& 1& 0& 0& 1& 0& 1& 1\\
0& 0& 0& 0& 1& 0& 0& 1& 0& 0& 0& 0& 1& 0& 0& 0& 1& 0& 0& 0& 0& 1& 1& 1\\
0& 0& 1& 0& 0& 0& 1& 0& 0& 1& 0& 0& 1& 0& 0& 0& 0& 0& 0& 1& 1& 0& 0& 1\\
0& 0& 0& 1& 0& 0& 0& 1& 1& 0& 0& 0& 0& 1& 0& 0& 0& 0& 1& 0& 0& 1& 0& 1\\
0& 1& 0& 0& 1& 0& 0& 0& 0& 0& 0& 1& 0& 0& 1& 0& 1& 0& 1& 0& 0& 1& 0& 0\\
0& 1& 0& 0& 0& 1& 0& 0& 0& 0& 1& 0& 0& 0& 0& 1& 0& 1& 0& 1& 1& 0& 0& 0\\
0& 0& 0& 0& 0& 0& 0& 0& 1& 1& 0& 0& 0& 0& 1& 1& 1& 1& 0& 0& 0& 0& 0& 1\\
0& 0& 0& 0& 0& 0& 0& 0& 0& 0& 1& 1& 0& 0& 0& 0& 1& 1& 1& 1& 0& 0& 1& 0
\end{array}\right]}$$
$\hspace{-2.5cm}A^{-1}=\displaystyle\frac{1}{21}{\scriptsize\left[\begin{array}{rrrrrrrrrrrrrrrrrrrrrrrr}
0&1&5&5&4&4&1&1&-1&-1&1&1&-1&-1&-4&-4&4&4&-2&-2&-4&-4&-6&2\\
1&0&5&5&-4&-4&-1&-1&-2&-2&-4&-4&1&1&1&1&4&4&-1&-1&4&4&2&-6\\
5&5&0&1&-4&1&-2&-1&1&-1&4&-4&-1&-2&4&-4&2&-6&1&-1&1&-4&4&4\\
5&5&1&0&1&-4&-1&-2&-1&1&-4&4&-2&-1&-4&4&-6&2&-1&1&-4&1&4&4\\
4&-4&-4&1&0&4&-1&4&4&1&-4&5&-1&-6&1&-4&-1&1&-2&2&5&1&-1&-2\\
4&-4&1&-4&4&0&4&-1&1&4&5&-4&-6&-1&-4&1&1&-1&2&-2&1&5&-1&-2\\
1&-1&-2&-1&-1&4&0&1&-4&1&2&1&5&4&-2&4&4&-4&5&-4&-1&-6&-4&1\\
1&-1&-1&-2&4&-1&1&0&1&-4&1&2&4&5&4&-2&-4&4&-4&5&-6&-1&-4&1\\
-1&-2&1&-1&4&1&-4&1&0&4&-1&-6&-4&5&4&-1&1&-4&1&5&2&-2&4&-4\\
-1&-2&-1&1&1&4&1&-4&4&0&-6&-1&5&-4&-1&4&-4&1&5&1&-2&2&4&-4\\
1&-4&4&-4&-4&5&2&1&-1&-6&0&1&4&-2&4&5&-2&-1&4&-1&-4&1&-1&1\\
1&-4&-4&4&5&-4&1&2&-6&-1&1&0&-2&4&5&4&-1&-2&-1&4&1&-4&-1&1\\
-1&1&-1&-2&-1&-6&5&4&-4&5&4&-2&0&1&1&2&-4&4&1&-4&-1&4&1&-4\\
-1&1&-2&-1&-6&-1&4&5&5&-4&-2&4&1&0&2&1&4&-4&-4&1&4&-1&1&-4\\
-4&1&4&-4&1&-4&-2&4&4&-1&4&5&1&2&0&1&-2&-1&-1&-6&5&-4&1&-1\\
-4&1&-4&4&-4&1&4&-2&-1&4&5&4&2&1&1&0&-1&-2&-6&-1&-4&5&1&-1\\
4&4&2&-6&-1&1&4&-4&1&-4&-2&-1&-4&4&-2&-1&0&1&1&-4&1&-1&5&5\\
4&4&-6&2&1&-1&-4&4&-4&1&-1&-2&4&-4&-1&-2&1&0&-4&1&-1&1&5&5\\
-2&-1&1&-1&-2&2&5&-4&1&5&4&-1&1&-4&-1&-6&1&-4&0&4&1&4&-4&4\\
-2&-1&-1&1&2&-2&-4&5&5&1&-1&4&-4&1&-6&-1&-4&1&4&0&4&1&-4&4\\
-4&4&1&-4&5&1&-1&-6&2&-2&-4&1&-1&4&5&-4&1&-1&1&4&0&4&-2&-1\\
-4&4&-4&1&1&5&-6&-1&-2&2&1&-4&4&-1&-4&5&-1&1&4&1&4&0&-2&-1\\
-6&2&4&4&-1&-1&-4&-4&4&4&-1&-1&1&1&1&1&5&5&-4&-4&-2&-2&0&1\\
2&-6&4&4&-2&-2&1&1&-4&-4&1&1&-4&-4&-1&-1&5&5&4&4&-1&-1&1&0
\end{array}\right]}$
\caption{The adjacency matrix and its inverse of a nuciferous graph with 24 vertices}\label{table2}
\end{table}
Based on our findings on nuciferous Cayley graphs we pose the following:
\begin{conj}
There exist infinitely many nuciferous Cayley graphs.
\end{conj}

\section*{Acknowledgments}
I would like to thank Ali Mohammadian for introducing the problem to me, and
anonymous referees for their comments which improved the presentation of the paper.
I'm also indebted to Patrick Fowler for double-checking the data given in Table~\ref{table} and pointing out an error in an earlier version of this paper. 
The research of the author was in part supported by a grant from IPM (No. 94050114).

{}

\end{document}